\documentclass[a4paper,14pt]{article}
\usepackage[utf8]{inputenc}
\usepackage{amssymb,amsmath,amscd}

\newcommand{\Sym}{\mathrm{Sym}}

\newcommand{\Sk}{\mathrm{Sk}}
\newcommand{\Q}{\mathbb{Q}}
\newcommand{\R}{\mathbb{R}}
\newcommand{\Z}{\mathbb{Z}}
\newcommand{\Cn}{\mathbb{C}}

\author{DMITRY V. GUGNIN\thanks{This work is supported by the Russian Science Foundation under grant 14-11-00414.}}

\title{\bf \large ON HOMEOMORPHISM TYPE OF SYMMETRIC PRODUCTS OF COMPACT RIEMANN SURFACES WITH PUNCTURES}

\date{}

\begin{document}

\pagestyle{plain}
\pagenumbering{arabic}

\maketitle

\begin{abstract}

Let $M^2_{g,k}$ and $M^2_{g',k'}$ be compact Riemann surfaces with punctures ($g,g'\ge 0$ --- genuses, $k,k'\ge 1$ --- number of punctures). For any Hausdorff space $X$ the quotient space $\Sym^nX :=  X^n/S_n$ is the $n$-th symmetric product of $X,  \   n\ge 2$. It is well known, that $\Sym^n M^2_{g,k}$ is a smooth quasi-projective variety. Open manifolds $\Sym^n M^2_{g,k}$ and $\Sym^n M^2_{g',k'}$ are homotopy equivalent {\it iff} $\  2g+k=2g'+k'$. 
\medskip

\textbf{Blagojevi\'{c}-Gruji\'{c}-\v{Z}ivaljevi\'{c} Conjecture (2003).} {\it Fix any $n\ge 2$, and two pairs $(g,k)$ and $(g',k')$ with the condition $2g+k=2g'+k'$. If $g\ne g'$, then open manifolds $\Sym^n M^2_{g,k}$ and $\Sym^n M^2_{g',k'}$ are not continuously homeomorphic.}
\medskip

The conjecture was proved in the paper \cite{Ziv03} by P.$\,$Blagojevi\'{c}, V.$\,$Gruji\'{c} and R.$\,$\v{Z}ivaljevi\'{c}   for the case $\mathrm{max}(g,g') \ge \frac{n}{2}$ (this implies the case $n=2$). As far as the author knows, up to this moment there were no results if $\mathrm{max}(g,g') < \frac{n}{2}$.
\medskip
 
The aim of this paper is to prove the conjecture in full generality.

\end{abstract}

\begin{flushleft}

\textbf{ \ \ \ \ \ \ \ \ \  2010 Mathematics Subject Classification.}  Primary 55S15; Secondary 57R19, 57R20. \\

\textbf{ \ \ \ \ \ \ \ \ \  Keywords.} Symmetric products, Riemann surface, characteristic classes, smoothings    \\ 
  \ \ \ \ \ \ \ \ \    \phantom{a}   of manifolds, attaching a boundary. 

\end{flushleft}

\medskip

\begin{flushright}
{\it \large To my son Alesha}
\end{flushright}

\section{Introduction}

Let $X$ be a Hausdorff space. The quotient space $\Sym^nX :=  X^n/S_n$ is called the $n$-th symmetric product of $X$ for any integer $n\ge 2$. It is easy to see that the functor $\Sym^n$ is a homotopy functor. If $X$ is a finite (countable) simplicial complex, then spaces $\Sym^nX, n \ge 2,$ admit at least two natural structures of a finite (countable) simplicial complexes. Therefore, if $X$ is homotopy equivalent to a countable (finite) CW-complex, then spaces $\Sym^n X, n \ge 2,$ are also homotopy equivalent to countable (finite) CW-complexes.\medskip

The (co)homology of $\Sym^n X$ was investigated by Nakaoka, Dold, Thom, Mattuck, Macdonald, Milgram and many other mathematicians starting from 1950-s. The main disadvantage of this functor is the following\medskip

\textbf{Fact $\alpha$.} {\it Suppose $M^m$ is any topological $m$-dimensional manifold, where $m\ge 3$. Then for any $n\ge 2$ the space $\Sym^n M^m$ is not even a homology manifold.}
\medskip

This fact is an easy consequence of the 
\medskip

\textbf{Fact $\beta$.} {\it For $m\ge 3$, the space $\Sym^2 \R^m$ is homeomorphic to $\R^m \times \mathrm{ConeInt}(\R P^{m-1})$. Here by $\mathrm{ConeInt}(X)$ we denote the open cone $\mathrm{Cone}(X) \backslash X$.}
\medskip

On the other hand, for 2-dimensional manifolds there symmetric products are always manifolds. This is the consequence of the  
\medskip

\textbf{Fact $\gamma$.} {\it There exists a canonical homeomorphism $\Sym^n \Cn \cong \Cn^n$ for all $n\ge 2$. \\
($n$ roots of a unital polynomial $\leftrightarrow$ its $n$ coefficients)} 
\medskip

This fact also implies that if we have an arbitrary Riemann surface $M^2$ (compact or noncompact), then the topological manifold $\Sym^n M^2, n\ge 2,$ inherits some natural structure of a complex $n$-dimensional manifold. Therefore, in this case we have a natural $C^{\omega}$ and $C^{\infty}$ structures on the manifold $\Sym^n M^2$ (these are the weakening of a holomorphic structure). 
\medskip

But, if we have two Riemann surfaces $M^2$ and $N^2$, that are only $C^{\infty}$-diffeomorphic, then smooth manifolds $\Sym^n M^2$ and $\Sym^n N^2$ {\it a priory} could not be $C^{\infty}$-diffeomorphic. We should recall the fundamental fact that 2-dimensional $C^{\infty}$-manifolds (closed or open) are $C^{\infty}$-diffeomorphic  {\it iff}  they are just  homeomorphic. Thus, we would like to pose the following 
\medskip

\textbf{Conjecture 1.} {\it Let $M^2$ and $N^2$ be Riemann surfaces, that are closed or closed with a finite number of punctures. Suppose that $M^2$ and $N^2$ are homeomorphic ($=$   $C^{\infty}$-diffeomorphic). Then smooth manifolds $\Sym^n M^2$ and $\Sym^n N^2$ are $C^{\infty}$-diffeomorphic for all $n\ge 2$. }
\medskip

The following fact is the main theorem of \cite{Mor}.
\medskip

\textbf{Fact $\delta$.} {\it The space $\Sym^n S^1, n \ge 2,$ is homeomorphic to the $D^{n-1}$-bundle over $S^1$, which is trivial for odd $n$, and non-oriented for even $n$. }
\medskip

We also need the following main result of \cite{O}.
\medskip

\textbf{Fact $\varepsilon$.} {\it The space $\Sym^n \bigvee\limits_{1}^{m} S^1, n \ge 2, m\ge 1,$ is homotopy equivalent to $\Sk^n T^m$. Here, the cell structure on the torus $T^m$ is the direct product of the minimal cell structure on $S^1$ (this structure has one $0$-cell).}
\medskip

The following fact is a folklore.
\medskip

\textbf{Fact $\zeta$.} {\it Suppose $X$ is a connected CW-complex. Then $\pi_1(\Sym^n X)= \pi_1(X)^{ab}=H_1(X;\Z)$ for all $n\ge 2$.} \medskip

Let us now focus on symmetric products of compact Riemann surfaces with a finite number of punctures. We denote by $M^2_{g,k}$ a Riemann surface of genus $g\ge 0$, and $k$ distinct points removed, $k\ge 1$.\medskip

It is obvious that $M^2_{g,k}$ is homotopy equivalent to the bouquet $\bigvee\limits_{1}^{s} S^1$, where $s=2g+k-1$. The fact\,$\zeta$ implies that $\pi_1(\Sym^n M^2_{g,k}) = \Z^{2g+k-1}$ for all $n\ge 2$. Let us fix any $n\ge 2$, and two pairs $(g,k)$ and $(g',k')$. From above, we have that open manifolds $\Sym^n M^2_{g,k}$ and $\Sym^n M^2_{g',k'}$ are homotopy equivalent {\it iff} $2g+k=2g'+k'$. 
\medskip

Now we are ready to formulate the
\medskip
 
\textbf{Blagojevi\'{c}-Gruji\'{c}-\v{Z}ivaljevi\'{c} Conjecture (2003).} {\it Fix any $n\ge 2$, and two pairs $(g,k)$ and $(g',k')$ with the condition $2g+k=2g'+k'$. If $g\ne g'$, then open manifolds $\Sym^n M^2_{g,k}$ and $\Sym^n M^2_{g',k'}$ are not continuously homeomorphic.}\medskip

This conjecture was posed by   Rade \v{Z}ivaljevi\'{c},  Pavle Blagojevi\'{c} and Vladimir Gruji\'{c} in \cite{Ziv03}, at  R.$\,$\v{Z}ivaljevi\'{c} lecture on the conference 18th British Topology Meeting, Manchester, September 2003, and also in \cite{Ziv05}.  According to \v{Z}ivaljevi\'{c}  (personal communication), they were very much inspired and influenced by the work of Kostadin Tren\v{c}evski and Don\v{c}o Dimovski \cite{TD1} and \cite{TD2}, see in particular the conjecture after Theorem A.2 in \cite{TD1}.  The conjecture was proved in \cite{Ziv03} for the case $\mathrm{max}(g,g') \ge \frac{n}{2}$ (this implies the case $n=2$). As far as the author knows, up to this moment there were no results if $\mathrm{max}(g,g') < \frac{n}{2}$. 
\medskip

The aim of this paper is to prove the following generalization of this conjecture.
\medskip

\textbf{Theorem 1.} {\it Fix any $n\ge 2$, and two pairs $(g,k)$ and $(g',k')$ with the condition $2g+k=2g'+k'$. If $g\ne g'$, then open manifolds $\Sym^n M^2_{g,k}\times\R^N$ and $\Sym^n M^2_{g',k'}\times\R^N$ are not continuously homeomorphic for all $N\ge 0$.}
\medskip

\textbf{Here is the plan of the proof of Theorem 1.}  
\medskip

Set $s:=2g+k-1=2g'+k'-1$. If $s = 0 \ \text{or} \ 1$, then $g=g'$. So, we have $s\ge 2$. Fix a pair $(g,k)$ with the condition $2g+k-1=s$. We will denote by $\Z_2$ the field $\Z/2\Z$.
\medskip

\textbf{Step 1.} The space $\Sym^n M^2_{g,k} \sim \Sk^n T^s$ has torsionless integral homology. Therefore, we have the ring isomorphism $H^*(\Sym^n M^2_{g,k};\Z_2) \cong H^*(\Sym^n M^2_{g,k};\Z)\otimes \Z_2$. 
\medskip

\textbf{Step 2.} The integral cohomology ring $H^*(\Sym^n M^2_{g,k};\Z)$ is equal to the cutted exterior algebra $\Lambda_{\Z}^{\le n}(\alpha_1,\alpha_2,\ldots,\alpha_s)$ for some $\Z$-basis $\alpha_1,\alpha_2,\ldots,\alpha_s$ of $H^1(\Sym^n M^2_{g,k};\Z)$. Thus, the ring $H^*(\Sym^n M^2_{g,k};\Z_2)$ is equal to $\Lambda_{\Z_2}^{\le n}(\overline{\alpha}_1,\overline{\alpha}_2,\ldots,\overline{\alpha}_s)$ for some (any) $\Z_2$-basis $\overline{\alpha}_1,\overline{\alpha}_2,\ldots,\overline{\alpha}_s$ of $H^1(\Sym^n M^2_{g,k};\Z_2)$.
\medskip

\textbf{Step 3.} The open manifold $\Sym^n M^2_{g,k}$ is a Zariski-open subset of the smooth projective variety $\Sym^n M^2_g$, where $M^2_g$ is the initial compact Riemann surface without punctures. The total Chern class of the complex manifold $\Sym^n M^2_g$ was computed by Macdonald in his famous paper \cite{Mac1}. The inclusion $i_{(n)}\colon \Sym^n M^2_{g,k} \to \Sym^n M^2_g$ induce the ring homomorphism $i^*_{(n)}\colon H^*(\Sym^n M^2_g;\Z) \to H^*(\Sym^n M^2_{g,k};\Z)$, and $i^*_{(n)}(c_1(\Sym^n M^2_g)) = c_1(\Sym^n M^2_{g,k})$.

From Macdonald's calculations one can easily derive that 
$$
c_1(\Sym^n M^2_{g,k}) = - (\alpha_1\smile\alpha_2 + \alpha_3\smile\alpha_4 + \ldots + \alpha_{2g-1}\smile\alpha_{2g})
$$ 
for some $\Z$-basis $\alpha_1,\alpha_2,\ldots,\alpha_s$ of $H^1(\Sym^n M^2_{g,k};\Z)$. 
\medskip

\textbf{Step 4.} Suppose we have a complex vector bundle $\xi\colon E\to B$ with the fiber $\Cn^n$ and the base $B$, which is a connected ENR (compact, or non-compact and homotopy equivalent to a finite polyhedron). Then the Stiefel-Whitney classes $w_k$ of the realization $\xi_{\R}\colon E_{\R}\to B$ can be computed from the Chern classes $c_l$ of the initial vector bundle $\xi\colon E\to B$ as follows:
$$
w_{2k+1}(\xi_{\R}) = 0 \ \ \forall k\ge 0;  \ \   w_{2k}(\xi_{\R}) = \rho_2 (c_k(\xi)) \ \ \forall k\ge 1.
$$
Here, $\rho_2\colon H^*(B;\Z) \to H^*(B;\Z_2)$ is the reduction homomorphism. 
\medskip

The statement of this step is a well known fact. 
\medskip

\textbf{Step 5.} Combining two previous steps, we have that
$$
w_2(\Sym^n M^2_{g,k})  =  \overline{\alpha}_1\smile\overline{\alpha}_2 + \overline{\alpha}_3\smile\overline{\alpha}_4 + \ldots + \overline{\alpha}_{2g-1}\smile\overline{\alpha}_{2g}
$$ 
for some $\Z_2$-basis $\overline{\alpha}_1,\overline{\alpha}_2,\ldots,\overline{\alpha}_s$ of $H^1(\Sym^n M^2_{g,k};\Z_2)$. As one has 
$$
H^2(\Sym^n M^2_{g,k};\Z_2) \cong \Lambda^2(H^1(\Sym^n M^2_{g,k};\Z_2)),
$$
we get 
$$
w_2(\Sym^n M^2_{g,k})  =  \overline{\alpha}_1\wedge\overline{\alpha}_2 + \overline{\alpha}_3\wedge\overline{\alpha}_4 + \ldots + \overline{\alpha}_{2g-1}\wedge\overline{\alpha}_{2g}.
$$ 
\medskip

\textbf{Step 6. (Topological invariance of Stiefel-Whitney classes for open smooth manifolds)} 

Suppose we have closed smooth connected manifolds $M^n$ and $N^n$. By  celebrated Wu formula, if $f\colon M^n \to N^n$ is a homotopy equivalence, then $f^*(w_k(N^n)) = w_k(M^n)$ for all $k\ge 1$. It is the famous {\it Homotopy invariance} of Stiefel-Whitney classes for closed manifolds. 
\medskip

But, the trivial example $M^2 = S^1\times \R^1$ and $N^2$ = (open M\"{o}bius strip) shows us that even $w_1$ is not a {\it homotopy} invariant for open manifolds. 
\medskip

Now we want to pose the following
\medskip 

\textbf{Conjecture 2.} {\it Suppose we have a purely continuous homeomorphism $f\colon M^n \to N^n$ of two open connected smooth manifolds, which are homotopy equivalent to a finite polyhedron. Then $f^*(w_k(N^n)) = w_k(M^n)$ for all $1\le k\le n$.}
\medskip

{\it Remark.} This conjecture is trivially true for $w_1$ (a loop preserve or change the orientation), and for $w_n=0$.  
\medskip

Below we will prove this conjecture for $w_2$ with the following additional condition: the abelian groups $H_1(M^n;\Z)$ and $H_2(M^n;\Z)$ are torsionless and $H_2(M^n;\Z)$ is generated by the images of continuous mappings of torus $T^2$ to $M^n$. 
\medskip

\textbf{Step 7.} Combining the steps 5 and 6, we get that the topological type of the open manifold $\Sym^n M^2_{g,k}$ determines the genus $g$. Moreover, as Stiefel-Whitney classes are invariant under taking the direct product with the euclidian spaces $\R^N, N\ge 0$, we conclude the proof of Theorem 1.

\section{Steps 1-2}

Suppose $Z$ is a finite connected CW-complex. Fix any $n\ge 2$ and the commutative ring $R$. Denote by $X$ the $n$-skeleton $\Sk^n(Z)$. It is evident that the inclusion $i\colon X\hookrightarrow Z$ induce the isomorphism $i^*\colon H^k(Z;R) \cong H^k(X;R)$ for all $0\le k \le n-1$. Suppose also that the algebraic boundary $\partial\sigma$ (with $\Z$ coefficients) of any $(n+1)$-dimensional cell $\sigma$ of $Z$ is equal to zero. Then, it is easy to see that $i^*\colon H^n(Z;R) \to H^n(X;R)$ is also an isomorphism.
\medskip

As the induced mapping $i^*\colon H^*(Z;R) \to H^*(X;R)$ is a ring homomorphism and $\dim X \le n$, we get that the ring $H^*(X;R)$ is just the $(n+1)$-cutted ring $H^{*\le n}(Z;R)$. Moreover, the mapping $i^*\colon H^*(Z;R) \to H^*(X;R)$ just cuts the $*\ge (n+1)$ part of $H^*(Z;R)$. 
\medskip

All the above requirements are satisfied for $Z=T^s$ (with standard minimal cell structure) and any $n\ge 2$. Therefore, we have proved the steps 1-2.

\section{Step 3}

Let $M^2_g$ be an arbitrary compact Riemann surface of genus $g$ and without punctures. It is well known that the ring $H^*(M^2_g;\Z)$ is torsionless. Also one can choose the $\Z$-basis $\gamma_1,\gamma_2,\ldots,\gamma_{2g}\in H^1(M^2_g;\Z)$ with the property 
$$
\gamma_i \gamma_j = 0 \ \text{unless} \ i-j = \pm g; \ \ \ \gamma_i\gamma_{i+g} = - \gamma_{i+g}\gamma_i = \delta = [M^2_g]  \  (1\le i\le g).
$$
\medskip

Macdonald in his famous paper \cite{Mac1} proved that the ring $H^*(\Sym^n M^2_g;\Z)$ is torsionless and gave an explicit description of this ring. But, Macdonald's proof contained several gaps. The full verification of Macdonald's theorem was made by the author in preprint \cite{Gug15} using only algebraic topology tools. Another verification, which uses heavily algebraic geometry, was made in 2002 by del Ba\~{n}o \cite{Ba}. 
\medskip

For $g=0$, one has $M^2_0 = \Cn P^1$ and $\Sym^n \Cn P^1 = \Cn P^n$. Also for any $k\ge 1$, we have $M^2_{0,k} = \Cn \backslash \{ \mu_1,\ldots, \mu_{k-1} \}$ and $\Sym^n M^2_{0,k}$ is an open domain in $\Cn^n$. Therefore, the open manifold $\Sym^n M^2_{0,k}$ is parallizable. So, suppose $g\ge 1$.
\medskip

One has the canonical projection $\pi_n\colon (M^2_g)^n \to \Sym^n M^2_g$, which induces the isomorphism 
$$
\pi_n^*\colon H^*(\Sym^n M^2_g;\Q) \cong H^*((M^2_g)^n;\Q)^{S_n}.
$$ 
Macdonald's theorem tells that the torsionless ring 
$$
H^*(\Sym^n M^2_g;\Z) \subset H^*((M^2_g)^n;\Z)^{S_n} = (H^*(M^2_g;\Z)^{\otimes n})^{S_n}
$$ 
has multiplicative generators 
$$
\xi_1:=\chi(\gamma_1), \xi_2:= \chi(\gamma_2),\ldots, \xi_g:= \chi(\gamma_g);  \  \xi'_1:=\chi(\gamma_{g+1}), \xi'_2:= \chi(\gamma_{g+2}),\ldots, \xi'_g:= \chi(\gamma_{2g})  \ \text{and} 
\  \eta:=\chi(\delta), 
$$
where
$$
\chi(\omega):= \omega\otimes 1\otimes \ldots \otimes 1 + 1\otimes \omega\otimes 1\otimes \ldots \otimes 1 + \ldots + 1\otimes \ldots \otimes 1\otimes \omega  \ \ \text{for all} \  \omega \in H^*(M^2_g;\Z). 
$$

Macdonald also computed the total Chern class of the complex manifold $\Sym^n M^2_g$ (see \cite{Mac1}, theorem (14.5)):
$$
c(\Sym^n M^2_g)  =  (1+\eta)^{n-2g+1}(1 + \eta - \xi_1\xi'_1)(1 + \eta - \xi_2\xi'_2) \ldots (1 + \eta - \xi_g\xi'_g).
$$ 
\medskip
The open manifold $\Sym^n M^2_{g,k}$ is obviously a Zariski-open subset of the smooth projective variety $\Sym^n M^2_g$. The inclusion $i_{(n)}\colon \Sym^n M^2_{g,k} \to \Sym^n M^2_g$ induce the ring homomorphism $i_{(n)}^*\colon H^*(\Sym^n M^2_g;\Z) \to H^*(\Sym^n M^2_{g,k};\Z)$, and $i_{(n)}^*(c_1(\Sym^n M^2_g)) = c_1(\Sym^n M^2_{g,k})$. 
\medskip

Let $i\colon M^2_{g,k} \to M^2_g$ be the inclusion mapping. Then $i_n:= (i)^n\colon (M^2_{g,k})^n \to (M^2_g)^n$ and $i_{(n)}:= \Sym^n i\colon \Sym^n M^2_{g,k} \to \Sym^n M^2_g$ are the corresponding mappings. 
\medskip

Let us made the following auxiliary observation. Suppose $X$ and $Y$ are connected ENR's (compact, or non-compact and homotopy equivalent to a finite polyhedron). Suppose also that homology $H_*(X;\Z), H_*(Y;\Z)$ and $H_*(\Sym^n X;\Z), H_*(\Sym^n Y;\Z)$ are torsion-free for some fixed $n\ge 2$.  Let $i\colon X \to Y$ be a continuous mapping. One has the following commutative diagram:
$$
\begin{CD}
X\times\ldots\times X  @>i_{n}>> Y\times \ldots\times Y  \\
@VVV   @VVV  \\
\Sym^n X @>i_{(n)}>> \Sym^n Y
\end{CD}
$$
Let us use the notation 
$$
\chi(\omega) :=  \omega\otimes 1\otimes \ldots \otimes 1 + 1\otimes \omega\otimes 1\otimes \ldots \otimes 1 + \ldots + 1\otimes \ldots \otimes 1\otimes \omega.  
$$
If $\omega\in H^*(X;\Z)$, then symmetric tensor $\chi_X(\omega) \in H^*(X^n;\Z)$ lies in the subring $H^*(\Sym^n X;\Z)$. (This fact follows from Nakaoka's Theorem (2.7) in \cite{Nak}. Also it was rediscovered by the author (see \cite{my5}, Integrality Lemma)).\\
For the above diagram one has $i^*_{(n)} (\chi_Y(\omega)) = \chi_X (i^*(\omega))$ for all $\omega\in H^*(Y;\Z)$. So, for the case $X :=  M^2_{g,k}$ and $Y := M^2_g$ we get 
$$
i^*_{(n)}(\eta)  = i^*_{(n)} (\chi(\delta)) = \chi (i^*(\delta)) = \chi(0) = 0.
$$
Therefore, one has the formula for the total Chern class of the manifold $\Sym^n M^2_{g,k}$:
$$
c(\Sym^n M^2_{g,k})  =  (1  - i^*_{(n)}(\xi_1)i^*_{(n)}(\xi'_1))(1  - i^*_{(n)}(\xi_2)i^*_{(n)}(\xi'_2)) \ldots (1 - i^*_{(n)}(\xi_g)i^*_{(n)}(\xi'_g)).
$$
This formula implies the presentation for the first Chern class
$$
c_1(\Sym^n M^2_{g,k})  \   =  \  - [ i^*_{(n)}(\xi_1)i^*_{(n)}(\xi'_1)  +  i^*_{(n)}(\xi_2)i^*_{(n)}(\xi'_2) +  \ldots +  i^*_{(n)}(\xi_g)i^*_{(n)}(\xi'_g) ] \   =  
$$
$$
-  \  [ \chi(i^*(\gamma_1))\chi(i^*(\gamma_{g+1})) + \chi(i^*(\gamma_2))\chi(i^*(\gamma_{g+2})) + \ldots + \chi(i^*(\gamma_{g}))\chi(i^*(\gamma_{2g})) ]. 
$$

Now we need a one more observation. Let $X$ and $Y$ be as above. Suppose that the induces homomorphism $i_*\colon \pi_1(X) \to \pi_1(Y)$ is an epimorphism. This implies that $i_*\colon H_1(X;\Z) \to  H_1(Y;\Z)$ is also an epimorphism. As the abelian groups $H_1(X;\Z)$ and $H_1(Y;\Z)$ are torsionless, one has that $i^*\colon  H^1(Y;\Z)  \to  H^1(X;\Z)$ is a monomorphism and the image $i^*(H^1(Y;\Z))$ is the direct summand of $H^1(X;\Z)$. 
\medskip

For our case $X :=  M^2_{g,k}$ and $Y := M^2_g$ the above condition is satisfied. So, we have that $i^*(\gamma_1), i^*(\gamma_2), \ldots, i^*(\gamma_{2g})$ is a part of some $\Z$-basis $\varepsilon_1, \ldots, \varepsilon_{2g}, \varepsilon_{2g+1}, \ldots, \varepsilon_{s}$ of the free abelian group $H^1(M^2_{g,k};\Z)$.  
\medskip

The theorem 1 from the author's paper \cite{Gug15} implies the following
\medskip

\textbf{Fact $\eta$.} {\it Let $X$ be a connected ENR, compact or non-compact and homotopy equivalent to a finite polyhedron. Suppose $\varepsilon_1, \varepsilon_2, \ldots, \varepsilon_{s}$ is a $\Z$-basis of the free abelian group $H^1(X;\Z)$. Then elements $\chi(\varepsilon_1), \ldots, \chi(\varepsilon_s)$ form a $\Z$-basis of the free abelian group  $H^1(\Sym^n X;\Z)$ for all $n\ge 2$.} 
\medskip

Due to this fact and the above observations we get that
$$
c_1(\Sym^n M^2_{g,k}) = - (\alpha_1\smile\alpha_2 + \alpha_3\smile\alpha_4 + \ldots + \alpha_{2g-1}\smile\alpha_{2g})
$$ 
for some $\Z$-basis $\alpha_1,\alpha_2,\ldots,\alpha_s$ of $H^1(\Sym^n M^2_{g,k};\Z)$. The Step 3 is proved.

\section{Step 5}

As the homology $H_*(\Sym^n M^2_{g,k};\Z)$ is torsionless, one has $H^*(\Sym^n M^2_{g,k};\Z_2) = H^*(\Sym^n M^2_{g,k};\Z) \otimes \Z_2$.
Combining two previous steps, we have that
$$
w_2(\Sym^n M^2_{g,k})  =  \overline{\alpha}_1\smile\overline{\alpha}_2 + \overline{\alpha}_3\smile\overline{\alpha}_4 + \ldots + \overline{\alpha}_{2g-1}\smile\overline{\alpha}_{2g}
$$ 
for some $\Z_2$-basis $\overline{\alpha}_1,\overline{\alpha}_2,\ldots,\overline{\alpha}_s$ of $H^1(\Sym^n M^2_{g,k};\Z_2)$. As one has 
$$
H^2(\Sym^n M^2_{g,k};\Z_2) \cong \Lambda^2(H^1(\Sym^n M^2_{g,k};\Z_2)),
$$
we get 
$$
w_2(\Sym^n M^2_{g,k})  =  \overline{\alpha}_1\wedge\overline{\alpha}_2 + \overline{\alpha}_3\wedge\overline{\alpha}_4 + \ldots + \overline{\alpha}_{2g-1}\wedge\overline{\alpha}_{2g}.
$$

\section{Step 6}

Above we posed the following
\medskip 

\textbf{Conjecture 2.} {\it Suppose we have a purely continuous homeomorphism $f\colon M^n \to N^n$ of two open connected smooth manifolds, which are homotopy equivalent to a finite polyhedron. Then $f^*(w_k(N^n)) = w_k(M^n)$ for all $1\le k\le n$.}
\medskip

{\it Remark.} This conjecture is trivially true for $w_1$ (a loop preserve or change the orientation), and for $w_n=0$.  
\medskip

Now we will prove this conjecture for $w_2$ with the following additional condition: \\
{\it the abelian groups $H_1(M^n;\Z)$ and $H_2(M^n;\Z)$ are torsionless and $H_2(M^n;\Z)$ is generated by the images of continuous mappings of torus $T^2$ to $M^n.$}
\medskip

As Stiefel-Whitney classes for smooth manifolds are trivially invariant under taking the direct product with euclidian spaces $\R^N, N\ge 0,$ we can assume that the dimension $n$ is as big as we want. Suppose that $n\ge 6$. 
\medskip

By the above condition on $H_1(M^n;\Z)$ and $H_2(M^n;\Z)$ one has that \\ 
(1) $H_2(M^n;\Z_2) = H_2(M^n;\Z)\otimes \Z_2$ and \\
(2) there exists a $\Z_2$-basis $\alpha_1 = g_{1*}[T^2],\alpha_2 = g_{2*}[T^2],\ldots,\alpha_t = g_{t*}[T^2]$ in $H_2(M^n;\Z_2)$ for some smooth embeddings $g_1,\ldots, g_t\colon T^2 \to M^n$. 
\medskip

Therefore, the topological invariance of $w_2$ in this case is a consequence of the following
 \medskip
 
\textbf{Lemma 1.} {\it Suppose we have a purely continuous homeomorphism $f\colon M^n \to N^n$ of two open connected smooth manifolds of dimension $n\ge 6$, which are homotopy equivalent to a finite polyhedron. Let $g\colon T^2 \to M^n$ be a smooth embedding. Then one has $w_{2M}(g_*[T^2]) = w_{2N}((fg)_*[T^2])$.} 
\medskip

\textbf{Proof.} We need the following auxiliary observation.
\medskip

Suppose $L^n$ is a connected smooth manifold, closed or open and homotopy equivalent to a compact polyhedron. Let $K^p$ be a connected smooth closed manifold of dimension $1\le p\le n-1$. Suppose $g\colon K^p \to L^n$ is a continuous mapping, and $g(K^p)\subset U^n$, where $U^n$ is an open domain in $L^n$ of finite homotopy type. Then $w_{pL} (g_*[K^p]) = w_{pU} (g_*[K^p])$, where $w_{pL}$ and $w_{pU}$ are the $p$-th Stiefel-Whitney classes of manifolds $L^n$ and $U^n$ respectively. This fact is well known. We will call this observation a Locality property. 
\medskip

Now we go back to the Lemma conditions. Let $U^n$ be a very small tubular neighbourhood of the submanifold $g(T^2)\subset M^n$. We take $U^n$ with the boundary $\partial U^n = L^{n-1}$, which is a smooth fiber bundle over $g(T^2)$ with the fiber $S^{n-3}$. One trivially has $\pi_1(L^{n-1}) = \pi_1(T^2) = \Z \oplus \Z$. 
\medskip

We get the continuous homeomorphism $f\colon U^n \to f(U^n) =: V^n\subset N^n$, where $\mathrm{int} (V^n)$ is an open domain in $N^n$. The open manifold $\mathrm{int} (V^n)$ inherits a smooth structure from $N^n$.
\medskip

From above observation we have that $w_{2M}(g_*[T^2]) = w_{2U}(g_*[T^2])$ and $w_{2N}((fg)_*[T^2]) = w_{2V}((fg)_*[T^2])$. Therefore, we need to prove the equality $w_{2U}(g_*[T^2]) = w_{2V}((fg)_*[T^2])$. Here $f\colon U^n \to V^n$ is a purely continuous homeomorphism, $U^n$ is a compact connected smooth manifold with the boundary, and $\mathrm{int} (V^n)$ has some smooth structure. Moreover, the boundary $\partial U^n = L^{n-1}$ is connected and has a free abelian fundamental group.  
\medskip

Let us take the doubles $\hat{U}^n:= U^n \cup_{L^{n-1}} U^n$ and $\hat{V}^n:= V^n \cup_{f(L^{n-1})} V^n$. The manifold $\hat{U}^n$ is smooth. The double $\hat{V}^n$ has two open domains with equal smooth structures (left part and right part), but a priory we have no natural smooth structure around the topologically locally flat codimension 1 submanifold $f(L^{n-1})\subset \hat{V}^n$.  
\medskip

The continuous homeomorphism $f\colon U^n \to V^n$ can be naturally extended to continuous homeomorphism $\hat{f}\colon \hat{U}^n \to \hat{V}^n$.
\medskip 

By the above Locality property and topological (even homotopy) invariance of Stiefel-Whitney classes for closed smooth manifolds, to conclude the proof of Lemma 1 it is sufficient to have the following
\medskip

\textbf{Lemma 2.} {\it Suppose $V^n$ is a topological connected compact $n$-manifold, $n\ge 6,$ with the connected boundary $\partial V^n$ such that $\pi_1(\partial V^n)$ is a free abelian group. Suppose also that $\Sigma$ is a smooth structure on the interior $\mathrm{int}(V^n)$. Take the double $\hat{V}^n:= V^n_{+}\cup_{\partial V^n} V^n_{-}$, where $V^n_{\pm}$ are the two copies of $V^n$. Let $\hat{\Sigma}$ be the union of the smooth structures $\Sigma_{\pm} = \Sigma$ on the union of open domains $\mathrm{int}(V^n_{+}) \cup \mathrm{int}(V^n_{-})$. Then for any $\varepsilon > 0$, there exists some smooth structure $\hat{\Sigma}_{\varepsilon}$ on the whole double $\hat{V}^n$ such that the structures $\hat{\Sigma}_{\varepsilon}$ and $\hat{\Sigma}$ coincides on the open subset $U_{\varepsilon}:= \{x\in \hat{V}^n | d(x, \partial V^n) > \varepsilon \} $. Here $d(\cdot,\cdot)$ is an arbitrary fixed metric on the metrizable compact $\hat{V}^n$.}
\medskip

\textbf{Proof.} One has a natural involution $\tau\colon \hat{V}^n \to \hat{V}^n$, which permutes the left and the right parts $V^n_{\pm}$, and $\tau(x) = x$ {\it iff} $x\in \partial V^n$. Let us denote by $L^{n-1}$ the boundary $\partial V^n$. 
\medskip

By the topological collaring theorem there exists a collar $\overline{C}_{+}\approx L^{n-1}\times [0,1) \subset V^n_{+}$. Let us consider the collar interior $C_{+}$, which is homeomorphic to $L^{n-1} \times (0,1)$. The smooth structure $\Sigma$ induce some smooth structure $\Theta_{+}$ on the open domain $C_{+}$. 
\medskip

By the celebrated Product Structure Theorem (see \cite{KS}, Essay I), there exist some smooth structure $\Theta_{0}$ on $L^{n-1}$ and a diffeomorphism $f_{+}\colon L^{n-1}_{\Theta_0} \times (0,1) \to (L^{n-1} \times (0,1))_{\Theta_{+}}$. By the involution $\tau$, we get the symmetric collar $\overline{C}_{-} \subset V^n_{-}$, the symmetric smooth structure $\Theta_{-}$ on the open domain $C_{-}$ and the symmetric diffeomorphism $f_{-}\colon L^{n-1}_{\Theta_0} \times (-1,0) \to (L^{n-1} \times (-1,0))_{\Theta_{-}}$.
\medskip

Fix any $\varepsilon > 0$. There exists sufficiently small $0< \delta < \frac{1}{2}$ such that if $x\in \hat{V}^n$ and $d(x,L^{n-1}) > \varepsilon$, then
$$
x\notin f_{+} (L^{n-1}_{\Theta_0} \times (0,\delta]) \sqcup L^{n-1} \sqcup f_{-} (L^{n-1}_{\Theta_0} \times [-\delta,0)).
$$
\medskip

Let us consider the subset $f_{+} (L^{n-1}_{\Theta_0} \times (0,\delta]) \sqcup L^{n-1}$. By standard argumentation it is a topological $h$-cobordism. But, the fundamental group $\pi=\pi_1(L^{n-1})$ is a free abelian group, so the Whitehead torsion $\mathrm{Wh}(\pi)$ is zero (it is a classical theorem of Bass-Heller-Swan). Thus, the $h$-cobordism $f_{+} (L^{n-1}_{\Theta_0} \times (0,\delta]) \sqcup L^{n-1}$ is an $s$-cobordism. By the $s$-cobordism theorem in topological category this $s$-cobordism is a cylinder. 
\medskip

Therefore, one has a homeomorphism $h_{+} \colon  L^{n-1}_{\Theta_0}\times\{ \delta \}  \to L^{n-1}=\partial V^n$ and a respective homeomorphism of cylinders 
$$
\overline{h}_{+} \colon f_{+} (L^{n-1}_{\Theta_0} \times (0,\delta]) \sqcup L^{n-1} \to  (L^{n-1}_{\Theta_0}\times\{ \delta \}) \times [0,1].
$$
\medskip

By the action of the involution $\tau$, we get the symmetric homeomorphism $h_{-} \colon  L^{n-1}_{\Theta_0}\times\{ -\delta \}  \to L^{n-1}=\partial V^n$ and a symmetric respective  homeomorphism of cylinders 
$$
\overline{h}_{-} \colon f_{-} (L^{n-1}_{\Theta_0} \times [-\delta,0)) \sqcup L^{n-1} \to  (L^{n-1}_{\Theta_0}\times\{ -\delta \}) \times [-1,0].
$$
\medskip

Now it is easy to see that we get the following decomposition of the manifold $\hat{V}^n$:
$$
\hat{V}^n = ( \mathrm{int}(V^n_{+})\backslash   f_{+} (L^{n-1}_{\Theta_0} \times (0,\delta)) ) \bigcup (L^{n-1}_{\Theta_0}\times [-1,1]) \bigcup  ( \mathrm{int}(V^n_{-})\backslash   f_{-} (L^{n-1}_{\Theta_0} \times (-\delta,0))). 
$$

On the left part of this decomposition there is the initial smooth structure $\Sigma_{+}$, on the right part --- the initial smooth structure $\Sigma_{-}$, and in the middle there is the smooth structure of the cylinder. Moreover, all these structures are compatible on the boundaries $L^{n-1}_{\Theta_0}\times \{ \pm 1\} $. Therefore, we get a smooth structure $\hat{\Sigma}_{\varepsilon}$ with the needed property. 
\medskip

Lemmas 2 and 1 are completely proved.   \ \ \ \ \ \ \ \ \ \ \ \ \ \ \ \ \ \ \ \ \ \ \ \ \ \ \ \ \ \ \ \ \ \ \ \ \ \ \ \ \ \ \ \ \ \ \ \ \ \ \ \ \ \ \ \ \ \ \ \ \ \ \ \ \ \ \ \ \ \ \ \ \ \ \ \  \ \  \  \  $\Box$
\medskip

The proof of the Step 7 is a standard exercise in linear algebra. Therefore, we conclude the proof of Theorem 1.

\section{Pontrjagin classes and attaching a boundary}

For a complex manifold $M^{2n}$, which is closed or open and of a finite homotopy type, there is a standard procedure to calculate its Pontrjagin classes from its Chern classes.  This calculation for $\Sym^nM^2_{g,k}$ is not hard and gives the following
\medskip

\textbf{Proposition 1.} {\it For any $n\ge 2, g\ge 0, k\ge 1,$ all integral Pontrjagin classes of the manifold $\Sym^nM^2_{g,k}$ are equal to zero.}
\medskip

The following fact is not hard to prove and is a folklore.
\medskip

\textbf{Fact $\theta$.} {\it Let us denote by $D^2$ the closed $2$-disk. Then for all $n\ge 2$ the space $\Sym^n D^2$ is continuously homeomorphic to the closed $2n$-disk $D^{2n}$. But, there is no natural smoothing of $\Sym^n D^2$.}   
\medskip

\textbf{Corollary 1.} {\it Suppose $\overline{M}^2$ is a compact $2$-manifold with the boundary. Then the space $\Sym^n \overline{M}^2$ is a compact $2n$-manifold with the boundary for all $n\ge 2$.} 
\medskip

This corollary implies that we can in TOP category naturally attach a boundary to the open manifold $\Sym^nM^2_{g,k}$. We just need to take the compact Riemann surface $\overline{M}^2_{g,k}$ with $k$ small open disjoined disks removed (the boundary of these disks can be taken of $C^{\omega}$ class). Then for all $n\ge 2$ we get the topological compact $2n$-manifold $\Sym^n\overline{M}^2_{g,k}$ such that its interior is just $\Sym^n(\mathrm{int}\overline{M}^2_{g,k})$. 
\medskip

Therefore, the interior $\mathrm{int} (\Sym^n\overline{M}^2_{g,k})$ has a natural structure of a complex manifold. But, there is no natural smooth structure on the whole $\partial$-manifold $\Sym^n\overline{M}^2_{g,k}$. Moreover, the following question naturally arises:
\medskip

\textbf{Question 1.} {\it Could the TOP $\partial$-manifold $\Sym^n\overline{M}^2_{g,k}$ be smoothable?} 
\medskip

It is obvious that the compact $2$-manifold $\overline{M}^2_{g,k}$ possesses a triangulation. For any compact polyhedron $K$ the space $\Sym^n K, n\ge 2,$ inherits some natural triangulation. So, the $\partial$-manifold $\Sym^n\overline{M}^2_{g,k}$ is a compact polyhedron. But, in dimension $4$ the well known fact states that a TOP $4$-manifold $L^4$  (with or without a boundary) is smoothable, if it is a compact polyhedron. Thus, for $n=2$ the above Question 1 has a positive answer. 
\medskip

\textbf{Proposition 2.} {\it For $n\ge 3$ the boundary $\partial \Sym^n\overline{M}^2_{g,k}$ is a connected closed manifold with a free abelian fundamental group.}
\medskip

\textbf{Proof.} Let us fix any $n\ge 3, g\ge 0, k \ge 1$. It is easy to show that the boundary $\partial \Sym^n\overline{M}^2_{g,k}$ is a connected closed manifold with an abelian fundamental group. The compact $\partial$-manifold $\Sym^n\overline{M}^2_{g,k}$ is orientable, so we can use the standard Poincar\'{e} Duality. 
\medskip

Let us use the following notation: $L^{2n}:= \Sym^n\overline{M}^2_{g,k}, \    K^{2n-1} :=  \partial \Sym^n\overline{M}^2_{g,k}$ and $\tilde{M}^{2n}:= L^{2n}/ K^{2n-1}$. 
\medskip

The Poincar\'{e} Duality with $\Z$-coefficients gives the commutative diagram: 
$$
\begin{CD}
H^{2n-2}(L^{2n})  @>>> H^{2n-2}(K^{2n-1}) @>>> H^{2n-1}(\tilde{M}^{2n})   @>>> H^{2n-1}(L^{2n})  \\
@VVV   @VVV  @VVV  @VVV  \\
H_2(\tilde{M}^{2n}) @>>> H_1(K^{2n-1}) @>>> H_1(L^{2n}) @>>> H_1(\tilde{M}^{2n}).
\end{CD}
$$
All vertical arrows are isomorphisms and rows are exact. 
\medskip

We know that $L^{2n}\sim \Sk^n T^s, \  s=2g+k-1$. So, for the dimension reasons we get $H^{2n-2}(L^{2n})  = 0$, and $H^{2n-1}(L^{2n}) = 0$. From the above diagram we also  have $H_2(\tilde{M}^{2n}) = H_1(\tilde{M}^{2n}) = 0$. Therefore, the inclusion $j\colon K^{2n-1} \hookrightarrow L^{2n}$ induces the isomorphism $j_*\colon H_1(K^{2n-1}) \cong H_1(L^{2n}) = \Z^{s}$. Above we mentioned that the fundamental group $\pi_1(K^{2n-1})$ is abelian. So, we get $\pi_1(K^{2n-1}) = H_1(K^{2n-1}) = \Z^s$. The proposition 2 is completely proved.     \ \ \ \ \ \ \ \ \ \ \ \ \ \ \ \ \ \ \ \ \ \ \ \ \ \ \ \ \ \ \ \ \ \ \ \ \ \ \ \ \ \ \ \ \ \ \ \ \ \ \ \ \ \ \ \ \ \ \ \ \ \ \ \ \ \ \ \ \ \ \ \ \ \ \ \  \ \  \  \  \  \ \ \ \ \ \ \ \ \ \ \ \ \ \ \ \ \ \ \ \ \ \ \ \ \ \ \ \ \ \ \ \ \ \ \  \   $\Box$
\medskip

Now we can answer the above Question 1 for any $n\ge 3$.
\medskip

\textbf{Proposition 3.} {\it The topological $\partial$-manifold $\Sym^n\overline{M}^2_{g,k}, n\ge 3,$ is smoothable. Let us denote by $\Sigma$ the natural smooth structure on the interior $\mathrm{int}(\Sym^n\overline{M}^2_{g,k})$. Then for any $\varepsilon > 0$ there exists a smooth structure $\hat{\Sigma}_{\varepsilon}$ on the whole manifold $\Sym^n\overline{M}^2_{g,k}$ such that this structure coincides with the structure $\Sigma$ on the open subset $U_{\varepsilon}:= \{x\in \Sym^n\overline{M}^2_{g,k}  \   |  \   d(x, \partial \Sym^n\overline{M}^2_{g,k}) > \varepsilon \} $. Here $d(\cdot,\cdot)$ is an arbitrary fixed metric on the metrizable compact $\Sym^n\overline{M}^2_{g,k}$.}
\medskip

To prove this Proposition one just has to take Proposition 2 and an evident version of Lemma 2 above.  $\Box$
\medskip

Let us fix any $n\ge 2,  \   g,g'\ge 0, \   k,k'\ge 1$, such that $2g+k = 2g'+k'$ and $g\ne g'$. At the end of this paper we want to pose the following
\medskip

\textbf{Conjecture 3.} {\it The smoothable topological manifolds $\partial \Sym^n\overline{M}^2_{g,k}$ and $\partial \Sym^n\overline{M}^2_{g',k'}$ are not homotopy equivalent.} 
\medskip

\centerline{ \large ACKNOWLEDGEMENTS}   

\bigskip

The author is deeply grateful to R.\v{Z}ivaljevi\'{c} for posing this Conjecture to the author. The author is also deeply grateful to his Advisor V.M.Buchstaber, to A.A.Gaifullin, T.E.Panov, A.V.Chernavsky and S.A.Melikhov for fruitful discussions.

\begin{flushleft}
{\it 
Steklov Mathematical Institute \\
of Russian Academy of Sciences, \\
Moscow, Russia \\
E-mail: dmitry-gugnin@yandex.ru
}
\end{flushleft}


\bigskip







\begin{thebibliography}{99}


\bibitem{Ziv03}	P.$\,$Blagojevi\'{c}, V.$\,$Gruji\'{c}, R.$\,$\v{Z}ivaljevi\'{c}, \textit{Symmetric products of surfaces and the cycle index}, Israel J. Math. \textbf{138} (2003), 61-72.

\bibitem{Ziv05} P.$\,$Blagojevi\'{c}, V.$\,$Gruji\'{c}, R.$\,$\v{Z}ivaljevi\'{c}, \textit{Arrangements of symmetric products of spaces}, Topology and its Applications \textbf{148} (2005), 213-232.


\bibitem{Ba} S.$\,$del$\,$Ba\~{n}o, \textit{On the Hodge theory of the symmetric powers of a curve}, Publ. Math. \textbf{46} (2002), 17-25.






\bibitem{my5}  D.$\,$V.$\,$Gugnin, \textit{Topological applications of graded Frobenius $n$-homomorphisms II}, Trudy Mos. matem. obschestva \textbf{73:2} (2012), 207-228 (in Russian). English translation: Trans. Moscow Math. Soc. (2012), 167-182.

\bibitem{Gug15} D.$\,$V.$\,$Gugnin, \textit{On integral cohomology ring of symmetric products}, (2015), arXiv preprint 1502.01862, 14 pp., submitted for publication in Moscow Math. Journal. 


\bibitem{KS} R.$\,$C.$\,$Kirby,$\,$L.$\,$C.$\,$Siebenmann, \textit{Foundational essays on topological manifolds, smoothings, and triangulations}, Princeton University Press, Princeton, New Jersey, 1977. 


\bibitem{Mac1}  I.$\,$G.$\,$Macdonald, \textit{Symmetric products of an algebraic curve}, Topology \textbf{1} (1962), 319-343.


\bibitem{Mac2}  I.$\,$G.$\,$Macdonald, \textit{The Poincar\'e polynomial of a symmetric product}, Proc. Camb. Phil. Soc. \textbf{58} (1962), 563-568.




\bibitem{Mor} H.$\,$R.$\,$Morton, \textit{Symmetric products of the circle}, Proc. Camb. Phil. Soc. \textbf{63} (1967), 349-352.


\bibitem{Nak} M.$\,$Nakaoka, \textit{Cohomology of symmetric products}, J. Inst. Polytech. Osaka City Univ. \textbf{8:2} Ser. A (1957), 121-145.



\bibitem{O} B.$\,$Ong, \textit{The homotopy type of the symmetric products of bouquets of circles}, International J. Math. \textbf{14:5} (2003), 489-497.




\bibitem{TD1} K.$\,$Tren\v{c}evski, D.$\,$Dimovski, \textit{On the affine and projective commutative $(m+k,m)$-groups}, J. of Algebra \textbf{240} (2001), 338-365.


\bibitem{TD2} K.$\,$Trenchevski$\,$and$\,$D.$\,$Dimovski, \textit{Complex Commutative Vector Valued Groups} (monograph), Macedonian Academy Sci. and Arts 1992.




\end{thebibliography}
\end{document}